\documentclass[10pt]{article}

\title{On the Significance of Digits in Interval Notation}
\author{\mbox{M.H. van Emden}\\
        \mbox{Computer Science Dept, University of Victoria}\\
        \mbox{Technical Report DCS-270-IR}
       }
\date{}

\newtheorem{definition}{Definition}

\begin{document}
\maketitle

\begin{abstract}
To analyse the significance of the digits used for interval bounds, we
clarify the philosophical presuppositions of various interval
notations. We use information theory to determine the information
content of the last digit of the numeral used to denote the interval's
bounds. This leads to the notion of \emph{efficiency} of a decimal digit:
the actual value as percentage of the maximal value of its information
content. By taking this efficiency into account, many presentations of
intervals can be made more readable at the expense of negligible loss of
information.
\end{abstract}

\section{Introduction}

Once upon a time, it was a matter of professional ethics among computers
never to write a meaningless decimal. Since then computers have become
machines and thereby lost any form of ethics, professional or
otherwise.  The human computers of yore were helped in their ethical
behaviour by the fact that it took effort to write spurious decimals.
Now the situation is reversed: the lazy way is to use the default
precision of the I/O library function. As a result, we are deluged with
meaningless decimals.

Of course interval arithmetic is not guilty of such negligence. After
all, the very {\em raison d'\^etre} of the subject is to be explicit
about the precision of computed results.  Yet, even interval arithmetic
is plagued by superfluous decimals, albeit in a more subtle way.  In
this note we first review the various interval notations. We argue in
favour of a rarely used notation called ``tail'', or ``factored'', which
has the advantage of avoiding the repetition of decimals that are
necessarily the same.  We analyse the information content of the
remaining decimals.

\section{Philosophical implications of an interval notation} 
Several papers \cite{hvnn01,szwc99,vnmdn01a} discussing interval
notations have been published recently.  The various notations have
different implications, just as people have different reasons for being
interested in interval arithmetic.

For some, intervals are a way of denoting a fuzzy, or perhaps
probabilistic, quantity. Others use intervals to give an indication
of the extent to which rounding has introduced error in a
computation. Here we assume an interpretation of intervals that does
not necessarily negate the above interpretations, but differs in the
way it is made precise. We call it the {\em set interpretation of
interval arithmetic}.

\paragraph{The set interpretation}
According to the set interpretation, variables range over the real
numbers. These reals are represented in computer memory as sets of
reals. The constraint is that if variable $x$ is represented by set
$S$, we have $x \in S$. Thus the set interpretation differs from
conventional numerical analysis in the {\em absence of errors}.  It is
either true or false that $x$ belongs to $S$.

The fact that $S$ contains more than one real is not an error.  In
conventional numerical analysis, an error arises when, for example, a
real variable $x$ with value $0.1$ is represented by a floating-point
number $f$. An error arises because $x = f$ is false. On the other
hand, representing $x$ by $S$ is not an error if $x \in S$.

Of course, the statement $x \in S$ provides only a limited amount
of information about $x$. The larger $S$ is, the less information. In
the set interpretation of interval arithmetic we distinguish
error, which is avoidable, from the inescapable fact that the
amount of information yielded by a finite machine is finite.

\paragraph{Consequences of the set interpretation}
Interval arithmetic is no exception to the rule that finite machines
can only give a finite amount of information.  In interval arithmetic
the sets of reals are limited to those that are easily representable:
closed, connected sets of reals that have finite floating-point numbers
as bounds, if they have a bound at all. Unbounded closed connected sets
of reals use the infinities of the floating-point standard in the
obvious way. Each of this finite set of sets of reals can be
represented by a pair of floating-point numbers. It is also the case
that for every set of reals, there exists a unique least floating-point
interval containing it.

This is the set interpretation of interval arithmetic. Its virtues
include that it is familiar. In fact, many people are surprised to
hear it given a name, as this is what they always thought intervals to
be. Another virtue is that, if the set interpretation is followed up
in all its consequences, it allows resolution of potential ambiguities
in interval arithmetic, especially in interval division involving
unbounded intervals, intervals containing zero, or intervals
containing nothing but zero \cite{hckvnmdn01}.

\section{Interval notations}
If one accepts the advantages of the set interpretation of interval
arithmetic, then one prefers a notation for an interval that suggests
a set. The traditional notation, exemplified by $[1.233,1.235]$ has this
advantage. Although widely used, it is not practical, as is
apparent\footnote{
I'm not making this up; see page 122 of \cite{vhlmyd97}.
} from the statement that an unknown real $x$ belongs to
\begin{eqnarray}
      [+0.6180339887498946804,+0.6180339887498950136]. \label{crude}
\end{eqnarray}
The problem with this ubiquitous notation is that it is hard to separate
two important pieces of information: \emph{where} the interval is, and
\emph{how wide} it is. To remedy this defect,
Hyv\"onen \cite{hvnn01} described a notation according to which one writes
instead
\begin{eqnarray}
      +0.61803398874989[46804,50136].                \label{factored}
\end{eqnarray}
The situation is similar when we are annoyed by having to write
\begin{eqnarray*}
      0.61803398874989x+0.61803398874989y,
\end{eqnarray*}
which we prefer to have in factored form: $0.61803398874989(x+y)$.
Hence we propose to refer to (\ref{factored}) as \emph{factored notation}
for intervals\footnote{
The notation has been occasionally used without comment in the literature;
see for example \cite{vnmdn99a}.
Credit goes to Hyv\"onen, whose paper \cite{hvnn01} was the first to
appear in print that drew attention to it and named it.
Independently I did so in \cite{vnmdn01a}.
Hyv\"onen called it ``tail notation''.
}.
The name is more than an analogy: in general,
one \emph{factors} with respect to a multiplicative infix operation,
of which concatenation on strings is an example.

In the example the bounds are in normalized scientific notation and
have the same exponent. In general, factored notation converts an
interval $[a \times 10^p,b \times 10^q]$, with normalized numerals as
bounds, first to $[a,b \times 10^{q-p}] \times 10^p$, where the upper
bound is not necessarily normalized.  When $p \not = q$, then this
cannot be shortened by taking an initial string of common first
decimals outside the brackets. It can only be shortened by limiting the
precision of $a$ and $b$, a topic we address later in the paper.

Table~\ref{CLASSIF} contains an overview of interval notations. Most
of the table is adopted from Hyv\"onen \cite{hvnn01}.
In this overview we distinguish three categories: (a) those that suggest
a set, (b) those that suggest a number degraded by an error, and (c) those
that suggest a pure number.
\begin{table}
\begin{tabular}{l|l|l|}
Notation  &          Interval value        &   Name of notation    \\
\hline
\hline
$[1.233,1.235]$ &    $[1.233,1.235]$       &      Classic          \\
$1.23[3,5]$     &    $[1.233,1.235]$       &      Factored         \\
$1.234\pm 2$    &    $[1.232,1.236]$       &      Range            \\
$1.234$\verb|~| &    $[1.2335,1.2345]$     &      Tilde            \\
$1.234+$        &    $[1.234,1.235]$       &      Plus             \\
$1.234+[-1e-3,2e-3]$&$[1.233,1.236]$       &      Error            \\
$1.234*$        &    $[1.233,1.235]$       &      Star             \\
$1.234$         &    $[1.233,1.235]$       &      Single-Number
\end{tabular}
\caption{
\label{CLASSIF}
Overview of interval notations, adapted from Hyv\"onen \cite{hvnn01}.
}
\end{table}
The Classic and Factored notations belong to category (a).
Under category (b) we have added, in analogy with the Tilde notation, the
Plus notation. This latter notation is useful in the improvement of the
factored notation discussed later on in this paper.
Category (c) is in the last line. Hyv\"onen used the
name ``Fortran notation''. The notation is actually the
``Single-number notation'' for the Fortran implementation described in
\cite{szwc99}.

The virtue of the notations in category (b) is that they make explicit
that a numeral is not to be interpreted according to mathematical notation,
by which we mean that
\begin{equation}
d_{m}d_{m-1}\ldots d_0.d_{-1}\ldots d_{-n}   \label{mathematical}
\end{equation}
denotes the number $\sum_{i= -n}^m d_i 10^i$.
Mathematical notation implies an infinite number of zeros after the last
digit when $n > 0$.

Mathematical notation is not the only way to interpret
(\ref{mathematical}).  For a long time physicists, chemists, and
engineers have used the convention that
\emph{a numeral has as meaning
any number that rounds to the number denoted by the numeral
displayed}.
The coexistence of mathematical notation with the physics convention
introduces an ambiguity that is often resolved by context.  With
intervals, the ambiguity becomes problematic, as we need numerals to
denote the bounds of an interval in the classic notation. Are these to
be interpreted according to mathematical notation, or according to the
physics convention?  It is implicit in most of the interval literature,
and explicit in \cite{hvnn01,vnmdn01a}, that \emph{the numerals in the
bounds of an interval are to be interpreted according to the
mathematical notation}. In this paper we follow that rule.

We therefore propose to avoid category (c) and to give single-number an
annotation to indicate that it does not have the usual mathematical
meaning. This has been done by Hickey, who introduced \cite{hck00} the
Star notation of Table~\ref{CLASSIF}.

\paragraph{Difficulties of factored notation}
There are two problems with the classical notation. The first is the
{\em scanning problem}\/: one needs to scan both bounds digit by digit
to find the leftmost different digit.  Only then does one have an idea
of the width of the interval. The second problem, the {\em problem of
useless digits} can also be found in (\ref{crude}): the width of the
interval is specified by no fewer than five digits.  Restricting
oneself to four digits for this purpose will give almost as much
information about $x$ and that the difference is so small as not to be
worth that fifth digit.  As we will show below, the same holds almost
always for all digits beyond the first two or three.

Factored notation solves the scanning problem; the problem of useless
digits remains. To solve it also, we need to study quantitatively the
information content of the statement that an unknown real $x$ is
contained in an interval $[a,b]$.

\section{Information theory}\label{INFOTH}
According to Shannon's theory of information (see for example, among
many textbooks, \cite{sh65}), observations can reduce the amount of
uncertainty about the value of an unknown quantity. The amount of
information yielded by an observation is \emph{the decrease (if any) in
the
amount of uncertainty}. Shannon argues that the amount of uncertainty
is appropriately measured by the \emph{entropy} of the probability
distribution over the possible values. For a uniform distribution on a
finite number of values, this reduces to the logarithm of the number
of possible values. It can be shown that the entropy for a
distribution over $n$ outcomes is maximized by the uniform
distribution over these outcomes.

When there are two equally probable possible values, and if one would
like this logarithm to come out at unity, one takes $2$ as base of the
logarithms and one calls the unit of information {\em bit}, for
{\bf b}inary un{\bf it} of information. Thus, the binary digits
carry at most one bit of
information. Similarly, if one works with decimal digits, then it is
convenient to use $10$ as the basis of the logarithms.
%By analogy, we
%call the resulting unit of information {\em dit},
%from {\bf d}ecimal un{\bf it} of information. For mathematical
%reasons, it is sometimes convenient to use $e \in 2.71[8,9]$ as basis
%for the logarithms, with the {\em nit} as corresponding unit of
%information.

Thus information theory determines for each number base the maximum
amount of information that can be carried by a digit. Normally, if we
don't know what a number is, and we are only given the first $k$
digits of a numeral denoting that number, we have no idea what the
next digit should be. That is, all possibilities in
$\{0,1,2,3,4,5,6,7,8,9\}$ are equally probable so that the uncertainty
is $\log_{10}10 = 1$. As a decimal digit can only distinguish between
ten possibilities, the efficiency of the $(k+1)$st digit is one.
%Of
%course, this also holds if we do know what the number is, and it is
%$\pi$, and if $k = 100$.

%Things are not as straightforward for the other important use of
%information theory in interval arithmetic. As we indicated before,
In the set interpretation of interval arithmetic, we have information of
the form that a real $x$ belongs to a set $S$. According to information
theory, this represents an uncertainty equal to the entropy of the
probability distribution over the elements of $S$.
What distribution to assume?
We are only interested in the large differences in
information carried by the successive digits of factored notation.
These are large compared to those due to the differences among
plausible distributions.

The fact that we are only interested in sets that are bounded
intervals, simplifies matters considerably. Plausible distributions for
bounded intervals include the uniform and the beta distributions.
From now on, if we know that $x$ is in
an interval $I$, we assume that the probability of $x$ belonging to any
subinterval of $I$ only depends on the width of that subinterval and
not on where in $I$ this subinterval is located. This property is
implied by the uniform distribution over $I$, and this is the
distribution we assume for computation of the uncertainty in the
statement $x \in I$. This uncertainty is equal to $- \log_{10} w$,
in decimal units of information, where $w$
is the width of $I$.

\section{Improvement of factored notation}

Factored notation solves the scanning problem. In this section we solve
the remaining problem that typically many of the digits inside the
brackets are useless. We do this by applying the formula found in
Section~\ref{INFOTH} to determine the information content of the digits
in factored notation. As factored notation is just an abbreviation of
it, this holds for classical notation as well.

We first consider a specific example in which we note a pattern of
rapidly decreasing efficiency as more digits are added. We explain this
phenomenon by a generally applicable formula, and use it to justify our
recommendation to write no more than three decimal digits inside the
brackets of factored notation.

For the example, we randomly selected an interval under the constraints
that both bounds have 15 digits, that the first five be the same, and
that the interval be nonempty. Thus we came to consider the interval
$[a,b]$ that is, in factored notation,
\begin{equation}
0.389015[282749894,960538227]            \label{SUPERF}
\end{equation}
The information content is $- \log_{10}(b-a)$, which is about 6.169
decimal units. If we have to represent the information that a real is
confined to this interval, but are only allowed to use two digits
inside the brackets, then this interval has to be $0.389015[28,97]$.
This interval has information content of about 6.161. Thus we saved
twice seven digits and lost an amount of information equal to 0.008
decimal units. Note that an optimally used pair of decimal digits in
factored notation carries 1.000 decimal units of information.

This example suggests that two decimals inside the brackets already
give almost all the information contained in the statement that $x$ is
in (\ref{SUPERF}).  That only two decimal digits inside the brackets
are enough could be a misleading feature of this particular example. To
investigate this possibility, we analyse the information content
remaining for all possible ways of shortening (\ref{SUPERF}).
From this we will see that a pattern emerges. We show that the pattern
is not a peculiarity of the example.
Because the pattern almost always occurs, we give it a name:
\emph{Rule of One Tenth}.
Before investigating this rule, we first need to be more
precise about shortening the representation of an interval.

\paragraph{Inflation}
Consider the statement that $x \in [a,b]$. Let $[a^\prime,b^\prime]$
properly contain $[a,b]$. Now it may be the case that $x \in
[a^\prime,b^\prime]$ conveys almost as much information about $x$ as $x
\in [a,b]$ and yet $[a^\prime,b^\prime]$ requires fewer digits to
write. Then $[a^\prime,b^\prime]$ is a more efficient representation
than $[a,b]$.

A more efficient representation such as $[a^\prime,b^\prime]$ 
may be obtained by one
or more applications of an operation we refer to as ``inflation''.
\begin{definition}
Let $I$ be the representation of an interval of which the bounds have a
finite number of decimals. The operation of \emph{inflation} has as
result the representation of the smallest interval containing $I$ where
each bound has one less decimal than the corresponding bound in $I$.
\end{definition}
In Table~\ref{EXAMPLES} we see some examples of inflation.
\begin{table}
\begin{tabular}{r||l|l}
 line number  & before inflation & after inflation              \\
\hline
\hline
      0       & $0.123[456,789]$  & $0.123[45,79]$                \\
      1       & $0.1[2345,34]$    & $0.1[234,4]$                  \\
      2       & $0.[1234,9999]$     & $[0.123,1.00]$                \\
      3       & $0.123[450,670]$  & $0.123[45,67]$                \\
      4       & $0.123[499,501]$  & $0.123[49,51]$
\end{tabular}
\caption{
Examples of inflation.
\label{EXAMPLES}
}
\end{table}
Line 0 is a typical case. Line 1 illustrates that inflation may apply
to intervals with an unequal number of decimals in the bounds. Line 2
is included to illustrate that inflation decreases the number of
digits, so that the four-digit $0.9999$ changes to the three-digit
numeral $1.00$.

Let us now consider the change in interval width due to inflation. In
line 3 of Table~\ref{EXAMPLES} we see that it can be as little as zero.
Line 4 shows that the width can increase by a factor of 10. In such a
case, the digits saved by inflation carry as much information as is
possible for a decimal digit.

In Table~\ref{SPREADSH} we see in the top line the bounds of interval
(\ref{SUPERF}).  Each next line shows the result of inflation applied
to the previous line. Thus it is true that $x$ is contained in each
interval of the table.
In the fourth column we see the information content of the
statement that $x$ belongs to the interval shown in that line. The last
column shows the decrease in information compared to the line before.
This decrease is to be compared to the information content of the
omitted decimal, which is $1$. Thus, the last column contains the
efficiency of showing the last decimal in each bound in the line
before.

\begin{table}
\begin{tabular}{r||l|l||l|l}
 &left boundary $a$ & right boundary $b$ & $-\log_{10}(b-a)$ & information \\
 &                &                    &           & loss                \\
\hline
\hline
0  & 0.389015 282749894 & 0.389015 960538227 & 6.168905911 &               \\
1  & 0.389015 28274989  & 0.389015 96053823  & 6.168905907 & 0.000000005   \\
2  & 0.389015 2827498   & 0.389015 9605383   & 6.168905804 & 0.000000103   \\
3  & 0.389015 282749    & 0.389015 960539    & 6.168904843 & 0.000000961   \\
4  & 0.389015 28274     & 0.389015 96054     & 6.168898435 & 0.000006407   \\
5  & 0.389015 2827      & 0.389015 9606      & 6.168834366 & 0.000064069   \\
6  & 0.389015 282       & 0.389015 961       & 6.168130226 & 0.000704140   \\
7  & 0.389015 28        & 0.389015 97        & 6.161150909 & 0.006979316   \\
8  & 0.38901 52         & 0.38901 60         & 6.096910013 & 0.064240896   \\
9  & 0.38901 5          & 0.38901 6          & 6           & 0.096910013   \\
10 & 0.3890 1           & 0.3890 2           & 5           & 1             \\
11 & 0.389 0            & 0.389 1            & 4           & 1             \\
12 & 0.3 89             & 0.3 90             & 3           & 1             \\
13 & 0.3 8              & 0.3 9              & 2           & 1             \\
14 & 0. 3               & 0. 4               & 1           & 1             \\
15 & 0                  & 1                  & 0           & 1          
\end{tabular}
\caption{
\label{SPREADSH}
Intervals $[a,b]$ containing an unknown real $x$.
Information loss as the result of successive inflations. Given that $x$
is in $[0,1]$, the information content of $x \in [a,b]$ is $-
\log_{10}(b-a)$. The loss due to inflation is in the last column.
}
\end{table}

As one goes down the table, considering successively more succinct, yet
true statements about $x$, one sees an interesting transition about
halfway. Of course something special has to happen at the point where
factored notation is $0.38901[5,6]$. The next more succinct intervals
are, successively, $0.3890[1,2]$, $0.389[0,1]$ and so on. In this
range, the information decrease is $1$, exactly the information content
of the decimal digit saved. That is, the digits that are saved here
are fully efficient.  Factored notation is not as useful here as it was
higher up in the table. In fact, it is redundant, as there is always a
pair of successive single decimals inside the brackets. An ad-hoc
notation in the style of tilde notation has a considerable advantage
here. I adopted the one proposed by Hickey \cite{hck00} and called
it ``Plus'' in Table~\ref{CLASSIF}.

Let us now consider the most important part of Table~\ref{SPREADSH}.
Suppose one considers shortening the interval in the top line to
$0.389015[28,97]$ and suppose one worries that too much information has
been lost. The last column in line 7 shows that the additional digits
contained in line 8 add only about one tenth of the amount information
contained in the last digits of line 7, which is already pretty low at
around one tenth of those in the line above that.  One can summarize
the last column above line 8 by the {\bf Rule of One Tenth}:
\begin{quote} \emph{Each additional digit carries about one tenth
of the information in the previous one.} \end{quote}
The rule holds quite well from line 8 upwards. If it would be exact,
the last column in line 1 would be $6*10^{-9}$ instead of the
$5*10^{-9}$ actually observed. Is this rule a fortuitous feature of
this particular example?  In the following, we will argue that it is
not.

\paragraph{The general case}
In Table~\ref{SPREADSH} we see that the Rule of One Tenth only holds
over many lines with considerable fluctuations from line to line. In
fact, in Table~\ref{EXAMPLES} we saw that inflation can cause an
increase in interval width of as little as a factor of one and as much
as a factor of ten.  These factors correspond to information losses of
0 and 1, respectively.  What can we say in general about interval
widening due to inflation?

\paragraph{}
We consider for the general case the interval shown digit by digit as
\begin{equation}\label{GENERIC}
0.x_1 \ldots x_{j-1}[y_j \ldots y_{j+k-1}p,z_j \ldots z_{j+k-1}q],
\end{equation}
where $y_j < z_j$ and $k \geq 2$.
We ask whether the number of digits can be safely decreased
by one application of the inflation operation.

If $p=q=0$, width does not increase, so inflation can be applied
without any loss of information. The largest information loss occurs if
$p=9$ and $q=1$, in which case the width increases by $18 \times
10^{-j-k}$.  Let us take $10^{-j-k+1}$ as a typical width increase, as
it is a convenient value near midway these extremes.

This increase should be compared with the width $w$ of
(\ref{GENERIC}).  The comparison is obscured by the large variation of
$w$. It may be as little as $10^{-j-k}$ (see last line of
Table~\ref{EXAMPLES}) and nearly as much as $10^{-j+1}$. In the case
(\ref{GENERIC}) is narrowest, inflation widens it typically by a factor
ten. In that case $p$ and $q$ carry as much information as is possible
for a decimal digit. Perhaps all decimals should be kept.  In the case
(\ref{GENERIC}) is widest, inflation widens it by a negligible amount.
Inflation is advisable.

Apparently it does not help to consider the extreme values of $w$, as
they lead to contradictory advice.  So let us consider average
values of $w$. We assume $k \geq 2$ (we retain at least two digits inside
the brackets).  If the average is in the order of $10^{-j}$, then
inflation causes negligible information loss.  If the average width is
near $10^{-j-k}$, then inflation causes the full amount of
information loss, so this is the worst case.  To simplify matters, we
make the worst case worse and assume that $w$ can range from $0$ to
$10^{-j+1}$. This is only a small change, as we are only interested in
$k > 2$, in which case the range from $0$ to $10^{-j-k}$ is negligible
compared to the range from $0$ to $10^{-j+1}$.

It is simplest to assume that the probability distribution of $w$ is
not far from uniform between $0$ and $10^{-j+1}$. In that case, it will
usually be the case that $w \in [10^{-j},10^{-j+1}]$.

%$10^{-j-k}$
%$10^{-j+1}$

But one may prefer not to make assumptions about the probability
distribution of $w$. Then one may accept the assumption that the digits
between the brackets in (\ref{GENERIC}) are independent random
variables with a uniform distribution on $\{0,\ldots,9\}$ under the
constraint that $y_j < z_j$. The average width of (\ref{GENERIC}) can
then be expressed as
\begin{equation}\label{AVERAGE}
w = \sum_{s = 0}^{9} \sum_{t = s+1}^{9} p_{st} w_{st}
\end{equation}
where $p_{st}$
is the probability of $y_j = s $ and $z_j = t$ and $w_{st}$ is the
average width under the constraint that $y_j = s$ and $z_j = t$.  For
$i$ between $0$ and $8$, if $y_j =i$, then $z_j$ can be $i,\ldots,9$.
Under the assumption about the distributions of the digits involved, we
have $p_{st} = 1/\sum_{i=1}^{9}i = 1/45$.

We are interested in a lower bound for $w_{st}$.
Each width is
bounded below by $(t-s-1)*10^{-j}$. Whatever the distribution, the
average is also bounded below by $(t-s-1)*10^{-j}$. Because this bound
depends only on $t-s$, we rewrite (\ref{AVERAGE}) as
\begin{eqnarray*}
w & = & \sum_{d=1}^{9} \sum_{a=0}^{9-d} p_{a,a+d} w_{a,a+d}
\end{eqnarray*}
Using $w_{a,a+d} \geq (d-1)*10^{-j}$ and $p_{st} = 1/45$ we have
\begin{eqnarray*}
w & \geq & (1/45) \sum_{d=1}^{9}(d-1)*10^{-j}  \\
  & \geq & (36/45)*10^{-j} = (4/5)*10^{-j}
\end{eqnarray*}
Moreover, $w$ is bounded above by $10^{-j+1}$. So it is
reasonable to assume that $w$ is in the order of $10^{-j}$.

%Let us now look at the effect of inflation acting on (\ref{GENERIC}).
%Dropping the
%last digit $p$ increases interval width by $p*10^{-j-k}$. Dropping $q$ and
%rounding upwards moves the upper bound by at most $(10-q)*10^{-j-k}$.
%Thus we have
%$$
%\begin{array}{rcccl}
%(10-9)*10^{-j-k} &\leq& (10-(q-p))*10^{-j-k} &\leq& (10 - (1-9))*10^{-j-k}
%								\\
%10^{-j-k}        &\leq& (10-(q-p))*10^{-j-k} &\leq& 18*10^{-j-k}
%\end{array}
%$$
%We take as typical value of interval widening $10^{-j-k+1}$, which is
%close to halfway these extremes.
Hence inflation widens an interval with a width of about $10^{-j}$ to
one that has a width of about $10^{-j}+10^{-j-k+1} =
10^{-j}(1+10^{-k+1})$.  Thus, the uncertainty decreased by the last
digit is in the order of $\log_{10}(1+10^{-k+1})$, which is about
$10^{-k+1}$, neglecting a factor of $\ln 10$.

This is also the decrease in information gain for every additional
digit inside the brackets in factored notation.  This is also the Rule
of Ten observed in Table~\ref{SPREADSH} when averaging over many rows.
We can expect that the third decimal in a factored notation only
increases information by $0.01$ of the potential information in a
decimal digit, and is therefore of questionable value. We recommend
factored notation with two decimals inside the brackets, while keeping
in mind that the rule does not apply in rare cases such as line $4$ in
Table~\ref{EXAMPLES}.

\section{Conclusions}

Interval methods are coming of age. When interval software was
experimental, it didn't matter whether interval output was easy to
read.  Now that the main technical challenges have been overcome, and
we at least \emph{know} how to ensure that the floating-point bounds
include all reals that are possible values of the variable concerned,
we need to turn our attention to small, mundane matters, which include
taking care of the convenience of users. Factored notation is an
advance in this respect.  However, without some attention to the number
of digits inside the brackets, one runs the risk of specifying in
maximum accuracy not the number under consideration, but the
unavoidable lack of information about this number.

\section{Acknowledgements}

I owe a debt of gratitude to the anonymous referees for their valuable
suggestions.
Many thanks to Fr\'ed\'eric Goualard for helpful comments on a draft of
this paper. We acknowledge generous support by the University of
Victoria, the Natural Science and Engineering Research Council NSERC,
the Centrum voor Wiskunde en Informatica CWI, and the Nederlandse
Organisatie voor Wetenschappelijk Onderzoek NWO.

\end{document}